\documentclass[12pt,oneside,reqno]{amsart}
\usepackage{txfonts}
\usepackage{bbm}
\usepackage{amsmath}
\usepackage{graphicx}
\usepackage{mathrsfs}
\usepackage{stmaryrd}
\usepackage{amsfonts}
\usepackage{enumerate,amsmath,amssymb,amsthm}
\newcommand{\be}{\begin{eqnarray}}
\newcommand{\ee}{\end{eqnarray}}
\newcommand{\ce}{\begin{eqnarray*}}
\newcommand{\de}{\end{eqnarray*}}
\newtheorem{theorem}{Theorem}[section]
\newtheorem{lemma}[theorem]{Lemma}
\newtheorem{remark}[theorem]{Remark}
\newtheorem{definition}[theorem]{Definition}
\newtheorem{proposition}[theorem]{Proposition}
\newtheorem{Examples}[theorem]{Example}
\newtheorem{corollary}[theorem]{Corollary}

\def\eps{\varepsilon}

\def\p{\partial}

\def\[{{\Big[}}
\def\]{{\Big]}}
\def\<{{\langle}}
\def\>{{\rangle}}
\def\({{\Big(}}
\def\){{\Big)}}

\def\bx{{\mathbf{x}}}

\def\dif{{\mathord{{\rm d}}}}

\def\={&\!\!=\!\!&}
\def\bt{\begin{theorem}}
\def\et{\end{theorem}}
\def\bl{\begin{lemma}}
\def\el{\end{lemma}}
\def\br{\begin{remark}}
\def\er{\end{remark}}
\def\bd{\begin{definition}}
\def\ed{\end{definition}}
\def\bp{\begin{proposition}}
\def\ep{\end{proposition}}
\def\bc{\begin{corollary}}
\def\ec{\end{corollary}}
\def\bx{\begin{Examples}}
\def\ex{\end{Examples}}

\def\cJ{{\mathcal J}}

\def\cL{{\mathcal L}}

\def\cS{{\mathcal S}}
\def\cT{{\mathcal T}}

\def\mB{{\mathbb B}}

\def\mE{{\mathbb E}}

\def\mH{{\mathbb H}}

\def\mN{{\mathbb N}}

\def\mR{{\mathbb R}}
\def\mS{{\mathbb S}}

\def\mW{{\mathbb W}}
\def\mX{{\mathbb X}}
\def\mY{{\mathbb Y}}

\def\sA{{\mathscr A}}
\def\sB{{\mathscr B}}

\def\geq{\geqslant}
\def\leq{\leqslant}

\def\iint{\int\!\!\!\int}

\allowdisplaybreaks

\begin{document}

\title{$L^p$-solvability of nonlocal parabolic equations with spatial dependent and non-smooth kernels$^*$}

\date{}
\author{Xicheng Zhang}


\thanks{$*$This work is supported by NSFs of China (No. 10971076) and Program
for New Century Excellent Talents in University (NCET 10-0654).}

\address{
School of Mathematics and Statistics,
Wuhan University, Wuhan, Hubei 430072, P.R.China\\
Email: XichengZhang@gmail.com
 }

\begin{abstract}
In this paper we prove the optimal $L^p$-solvability of nonlocal parabolic equation with spatial dependent and
non-smooth kernels.
\end{abstract}

\maketitle \rm

\section{Introduction}

In this paper we are considering the $L^p$-estimate of the following nonlocal operator:
\begin{align}
\cL^a f=\int_{\mR^d}[f(x+y)-f(x)-y^{(\alpha)}\cdot\nabla f(x)]a(x,y)|y|^{-d-\alpha}\dif y,\label{Lap}
\end{align}
where $\alpha\in(0,2)$, $a:\mR^d\times\mR^d\to\mR^+$ is a measurable function and
$$
y^{(\alpha)}:=1_{\alpha\in(1,2)}y+1_{\alpha=1}y1_{|y|\leq 1}.
$$
When $a(x,y)$ is {\it smooth} and {\it $0$-homogenous} in $y$, or $a(x,y)=a(y)$ is {\it independent} of $x$,
the $L^p$-estimates for this type of operators have been studied by Mikulevicius-Pragarauskas \cite{Mi-Pr1}
and Dong-Kim \cite{Do-Ki} (see also \cite{Zh1}).
However, for nonlinear applications, the smoothness and spatial-independence assumptions are usually not satisfied.

Let us now look at a nonlinear example. Consider the following variational integral appeared in nonlocal image and signal processing \cite{Gi-Os}:
$$
V(\theta):=\int_{\mR^d}\!\!\int_{\mR^d}\phi(\theta(x)-\theta(y))\kappa(x-y)|y-x|^{-d-\alpha}\dif x\dif y,\ \alpha\in(0,2),
$$
where $\phi:\mR\to\mR^+$ is an even convex $C^2$-function and $\kappa(-x)=\kappa(x)$. Assume that $\phi$ and $\kappa$
satisfy that for some $\Lambda>0$,
$$
\phi(0)=0,\ \ \ \Lambda^{-1}\leq\phi''(x)\leq\Lambda,
$$
and
$$
\Lambda^{-1}\leq \kappa(x)\leq\Lambda.
$$
The Euler-Lagrange equation corresponding to $V(\theta)$ is given by
$$
\int_{\mR^d}\phi'(\theta(t,y)-\theta(t,x))\kappa(y-x)|y-x|^{-d-\alpha}\dif y=0.
$$
In \cite{Ca-Ch-Va}, Caffarelli, Chan and Vasseur firstly considered the following time dependence problem:
$$
\p_t\theta(t,x)=\int_{\mR^d}\phi'(\theta(y)-\theta(x))\kappa(y-x)|y-x|^{-d-\alpha}\dif y,
$$
and proved that for any $\theta_0\in\mH^{1,2}$, there exists a unique global classical $C^{1,\beta}$-solution
to the above equation with $\theta(0,\cdot)=\theta_0$ in the $L^2$-sense. The existence of weak solutions
with non-increasing energy can be deduced by the standard energy argument. To address the regularity problem,
they followed the classical idea of De Giorgi and considered the following linearized equation
\begin{align}
\p_tw(t,x)=\int_{\mR^d}\phi''(\theta(t,y)-\theta(t,x))(w(t,y)-w(t,x))\kappa(y-x)|y-x|^{-d-\alpha}\dif y,\label{Eq1}
\end{align}
where $w(t,x)=\nabla\theta(t,x)$. If we set
$$
\hat k(t,x,y)=\phi''(\theta(t,y)-\theta(t,x))\kappa(y-x)
=\phi''\left((y-x)\cdot\int^1_0w(t,x+s(y-x))\dif s\right)\kappa(y-x),
$$
then, since $\phi''$ is an even function, we have
$$
\hat k(t,x,y)=\hat k(t,y,x),
$$
and equation (\ref{Eq1}) is understood in the weak sense: for all $\eta\in C^\infty_0(\mR^d)$,
\begin{align*}
\int_{\mR^d}\p_tw(t,x)\eta(x)\dif x=\int_{\mR^d}\!\!\int_{\mR^d}
(w(t,y)-w(t,x))(\eta(y)-\eta(x))\kappa(t,x,y)|y-x|^{-d-\alpha}\dif y\dif x.
\end{align*}
Clearly, if we let
$$
a(t,x,y):=\hat k(t,x,x+y),
$$
then equation (\ref{Eq1}) becomes
$$
\p_tw(t,x)=\int_{\mR^d}(w(t,x+y)-w(t,x))a(t,x,y)|y|^{-d-\alpha}\dif y.
$$
Notice that $a(t,x,y)$ is usually not smooth apriori in $x$ and $y$. This type of equation is our main motivation.

This paper is organized as follows: In Section 2, we give some necessary spaces. In Section 3, we
prove some estimates of nonlocal integral operators. In Section 4, the linear nonlocal parabolic
equation is studied. In a forthcoming paper, we shall use the result obtained in this paper to
study the stochastic differential equations with spatial dependence jump-diffusion coefficients (cf. \cite{Zh2}).

\textsc{Convention:} Throughout this paper, we shall use $C$ with or without subscripts to denote
an unimportant constant.

\section{Preliminaries}

In this section we introduce some necessary spaces of Dini-type (cf. \cite[p.30, (25)]{St}).
Let $\sA_0$ be the space of all  real bounded measurable functions $a:\mR^d\times\mR^d\to\mR$
with finite norm
$$
\|a\|_{\sA_0}:=\sup_{x,y\in\mR^d}|a(x,y)|+\int^1_0{\omega^{(0)}_a(r)\over r}\dif r<+\infty,
$$
where
\begin{align}
\omega^{(0)}_a(r):=\sup_{x\in\mR^d}\sup_{|y|\leq r}|a(x,y)-a(x,0)|.\label{Om2}
\end{align}
Let $\sA_1\subset\sA_0$ be the subspace with finite norm
$$
\|a\|_{\sA_1}:=\|a\|_{\sA_0}+\int^1_0{\omega^{(1)}_a(r)\over r}\dif r<+\infty,
$$
where
\begin{align}
\omega^{(1)}_a(r):=\sup_{|x-x'|\leq r}|a(x,0)-a(x',0)|.\label{Om22}
\end{align}

Let $\mN_0:=\mN\cup\{0\}$. For $p>1$ and $\beta\geq 0$, let $\mH^{\beta,p}:=(I-\Delta)^{-\frac{\beta}{2}}(L^p)$
be the Bessel potential space with the norm
$$
\|f\|_{\mH^{\beta,p}}:=\|(I-\Delta)^{\frac{\beta}{2}}f\|_p\sim\|f\|_p+\|(-\Delta)^{\frac{\beta}{2}}f\|_p,
$$
and for $q\in[1,\infty]$, let $\mB^{\beta,p}_q$ be the Besov space defined by
$$
\mB^{\beta,p}_q:=(L^p,\mH^{k,p})_{\frac{\beta}{k},q},
$$
where $k\in\mN$ and $\beta<k$, and $(\cdot,\cdot)_{\frac{\beta}{k},p}$ stands for the real interpolation space. Let us write
$$
\mW^{\beta,p}:=\mB^{\beta,p}_p.
$$
It is well-known that if $\beta$ is an integer and $p>1$, an equivalent norm in $\mW^{\beta,p}=\mH^{\beta,p}$ is given by
$$
\|f\|_{\mW^{\beta,p}}:=\sum_{k=0}^\beta\|\nabla^k f\|_p,
$$
where $\nabla^k$ denotes the $k$-order generalized gradient;
and if $0<\beta\not=\mbox{integer}$ and $p>1$, an equivalent norm in $\mW^{\beta,p}$ is given by
\begin{align}
\|f\|_{\mW^{\beta,p}}:=\|f\|_p+\sum_{k=0}^{[\beta]}
\left(\iint_{\mR^d\times\mR^d}\frac{|\nabla^kf(x)-\nabla^kf(y)|^p}
{|x-y|^{d+\{\beta\}p}}\dif x\dif y\right)^{\frac{1}{p}},\label{Eq404}
\end{align}
where for a number $\beta>0$, $[\beta]$ denotes the integer part of $\beta$
and $\{\beta\}:=\beta-[\beta]$. It is also well-known that Riesz's transform $\nabla(-\Delta)^{-\frac{1}{2}}$
is a bounded linear operator in $L^p$-space for any $p>1$ (see \cite{St}). Moreover, the following
interpolation inequality holds: for any $\beta\in(0,\gamma), p>1$ and $f\in\mH^{\gamma,p}$,
\begin{align}
\|(-\Delta)^{\frac{\beta}{2}}f\|_p\leq C\|f\|^{1-\frac{\beta}{\gamma}}_p
\|(-\Delta)^{\frac{\gamma}{2}}f\|^{\frac{\beta}{\gamma}}_p.
\end{align}
The following lemma is an easy consequence of \cite[Lemma 2.1]{Ko}.
\bl
For any $\beta\in(0,1)$, there exits a constant $C=C(\beta,d)>0$
such that for all $p\geq 1$ and $f\in\mH^{\beta,p}$,
\begin{align}
\|f(\cdot+y)-f(\cdot)\|_p\leq C|y|^\beta\|(-\Delta)^{\frac{\beta}{2}}f\|_p.\label{Le2}
\end{align}
\el

For each $t\in[0,1]$, write $\mY^{\beta,p}_t:=L^p([0,t];\mH^{\beta,p})$ with the norm
$$
\|u\|_{\mY^{\beta,p}_t}:=\left(\int^t_0\|u(s)\|^p_{\mH^{\beta,p}}\dif s\right)^{\frac{1}{p}},
$$
and let $\mX^{\beta,p}_t$ be the completion of all functions $u\in C^\infty([0,t];\cS(\mR^d))$ with
respect to the norm
$$
\|u\|_{\mX^{\beta,p}_t}:=\sup_{s\in[0,t]}\|u(s)\|_{\mH^{\beta-1,p}}+\|u\|_{\mY^{\beta,p}_t}+\|\p_tu\|_{\mY^{\beta-1,p}_t}.
$$
It is well-known that (cf. \cite[p.180, Theorem III 4.10.2]{Am}),
\begin{align}
\mX^{\beta,p}_t\hookrightarrow C([0,t];\mW^{\beta-\frac{1}{p},p}).\label{Em}
\end{align}
For simplicity of notation, we also write
$$
\mX^{\beta,p}:=\mX^{\beta,p}_1,\ \ \mY^{\beta,p}:=\mY^{\beta,p}_1.
$$

\section{$L^p$-estimate of nonlocal operators}

Let $\nu$ be a $\sigma$-finite measure on $\mR^d$, which is called a L\'evy measure if
$\nu(\{0\})=0$ and
$$
\int_{\mR^d}1\wedge|x|^2\nu(\dif x)<+\infty.
$$
Let $\Sigma$ be a finite measure on the unit sphere $\mS^{d-1}$ in $\mR^d$.
For $\alpha\in(0,2)$, define
\begin{align}
\nu^{(\alpha)}(B):=\int_{\mS^{d-1}}\left(\int^\infty_0
\frac{1_B(r\theta)\dif r}{r^{1+\alpha}}\right)\Sigma(\dif\theta),\ \ B\in\sB(\mR^d).\label{Eq4}
\end{align}
Then $\nu^{(\alpha)}$ is the L\'evy measure corresponding to the $\alpha$-stable process.
\bd
(i) Let $\nu_1$ and $\nu_2$ be two Borel measures on $\mR^d$. We say that $\nu_1$ is less than $\nu_2$ if
$$
\nu_1(B)\leq \nu_2(B),\ \ B\in\sB(\mR^d),
$$
and we simply write $\nu_1\leq \nu_2$ in this case.

(ii) The L\'evy measure $\nu^{(\alpha)}$ defined by (\ref{Eq4}) is called nondegenerate if
\begin{align}
\int_{\mS^{d-1}}|\theta_0\cdot\theta|^\alpha\Sigma(\dif\theta)\not=0,\ \ \forall\theta_0\in\mS^{d-1}.\label{Spe}
\end{align}
\ed

Throughout this paper we make the following assumption:
\begin{enumerate}[{\bf (H$^{(\alpha)}_\nu$)}]
\item Let $\nu$ be a L\'evy measure and satisfy that for some $\alpha\in(0,2)$,
\begin{align}
\nu^{(\alpha)}_1\leq \nu\leq\nu^{(\alpha)}_2,\ \ 1_{\alpha=1}\int_{r<|x|<R}y\nu(\dif y)=0,\ \
0<r<R<+\infty,\label{Ep0}
\end{align}
where $\nu^{(\alpha)}_i, i=1,2$ are two L\'evy measures with the form (\ref{Eq4}),
and $\nu^{(\alpha)}_1$ is nondegnerate.
\end{enumerate}

Let us recall the following result from \cite[Corollary 4.4]{Zh1}.
\bt\label{Th3}
Assume {\bf (H$^{(\alpha)}_\nu$)} with $\alpha\in(0,2)$. Then for any $p\in(1,\infty)$, there exists a constant
$C_0\in(0,1)$ such that for all $f\in\mH^{\alpha,p}$,
\begin{align}
C_0\|(-\Delta)^{\frac{\alpha}{2}}f\|_p\leq\|\cL^\nu f\|_p\leq C_0^{-1}\|(-\Delta)^{\frac{\alpha}{2}}f\|_p.\label{EW88}
\end{align}
\et

Below, for simplicity of notation, we write
\begin{align}
\cJ^{(\alpha)}_f(x,y):=f(x+y)-f(x)-y^{(\alpha)}\cdot\nabla f(x).\label{JJ}
\end{align}
We first prepare the following lemma for later use.
\bl\label{Le1}
Suppose that $a\in\sA_0$ and $\nu\leq\nu^{(\alpha)}$ for some $\alpha\in(0,2)$. For any $p>1$,
there exists a constant $C=C(\alpha,p,d)>0$ such that for all $f\in \mH^{\alpha,p}$ and $\eps\in(0,1)$,
$$
\left\|\int_{|y|\leq\eps}\cJ^{(\alpha)}_f(\cdot,y)
(a(\cdot,y)-a(\cdot,0))\nu(\dif y)\right\|_p\leq
C\|(-\Delta)^{\frac{\alpha}{2}}f\|_p\int^\eps_0\frac{\omega^{(0)}_a(r)}{r}\dif r,
$$
where $\omega^{(0)}_a$ is defined by (\ref{Om2}).
\el
\begin{proof}
Let us look at the case of $\alpha\in[1,2)$. Since $a\in\sA_0$, by Minkowski's inequality we have
\begin{align*}
&\left\|\int_{|y|\leq\eps}\Big[f(\cdot+y)-f(\cdot)-y\cdot\nabla f(\cdot)\Big]
(a(\cdot,y)-a(\cdot,0))\nu(\dif y)\right\|_p\\
&\quad\leq\int_{|y|\leq\eps}|y|\left(\int^1_0\|\nabla f(\cdot+sy)-\nabla f(\cdot)\|_p\dif s\right)
\omega^{(0)}_a(|y|)\nu^{(\alpha)}(\dif y)\\
&\quad\stackrel{(\ref{Le2})}{\leq} C\|(-\Delta)^{\frac{\alpha-1}{2}}\nabla f\|_p\int_{|y|\leq\eps}
|y|^\alpha\omega^{(0)}_a(|y|)\nu^{(\alpha)}(\dif y)\\
&\quad\leq C\|(-\Delta)^{\frac{\alpha}{2}} f\|_p\int^\eps_0\frac{\omega^{(0)}_a(r)}{r}\dif r,
\end{align*}
where the last step is due to (\ref{Eq4}) and the boundedness of Riesz transform in $L^p$-space.
The case of $\alpha\in(0,1)$ is similar.
\end{proof}

For $a\in\sA_0$, define the following nonlocal operator:
$$
\cL^{a\nu}f(x):=\int_{\mR^d}\cJ^{(\alpha)}_f(x,y)a(x,y)\nu(\dif y),
$$
where $\cJ^{(\alpha)}_f(x,y)$ is given by (\ref{JJ}).
We now establish the following characterization about the domain of $\cL^{a\nu}$.
\bt
Let $\alpha\in(0,2)$. Assume that {\bf (H$^{(\alpha)}_\nu$)} holds and
$a\in\sA_0$ satisfies that for some $0<a_0<a_1$ and any $0<r<R<\infty$,
\begin{align}
a_0\leq a(x,0)\leq a_1,\ \
1_{\alpha=1}\int_{r\leq |y|\leq R}ya(x,y)\nu(\dif y)=0.\label{EL22}
\end{align}
Then for any $p\in(1,\infty)$, there exists a constant $C_1\in(0,1)$
depending only on $a_0,a_1,\nu^{(\alpha)}_1,\nu^{(\alpha)}_2$ and $\alpha,d,p$ such that
for all $f\in\mH^{\alpha,p}$,
\begin{align}
C_1\|f\|_{\alpha,p}\leq\|\cL^{a\nu}f\|_p+\|f\|_p\leq C_1^{-1}\|f\|_{\alpha,p}.\label{EW1}
\end{align}
\et
\begin{proof}
We make the following decomposition:
\begin{align*}
\cL^{a\nu}f(x)&=a(x,0)\cL^\nu f(x)+\int_{|y|>\eps}\cJ^{(\alpha)}_f(x,y)(a(x,y)-a(x,0))\nu(\dif y)\\
&\quad+\int_{|y|\leq\eps}\cJ^{(\alpha)}_f(x,y)(a(x,y)-a(x,0))\nu(\dif y)\\
&=:I_1(x)+I_2(x)+I_3(x).
\end{align*}
For $I_1(x)$, by Theorem \ref{Th3} and condition (\ref{EL22}), we have
$$
a_0C_0\|(-\Delta)^{\alpha/2}f\|_p\leq\|I_1\|_p\leq a_1C_0^{-1}\|(-\Delta)^{\alpha/2}f\|_p.
$$
For $I_2(x)$,  if $\alpha=1$, by (\ref{EL22}) we have
\begin{align*}
\|I_2\|_p=\left\|\int_{|y|>1}[f(\cdot+y)-f(\cdot)](a(\cdot,y)-a(\cdot,0))\nu(\dif y)\right\|_p
\leq 4\|f\|_p\|a\|_\infty\nu(B^c_1);
\end{align*}
if $\alpha\in(0,1)$, we have
\begin{align*}
\|I_2\|_p\leq 4\|f\|_p\|a\|_\infty\nu(B^c_\eps);
\end{align*}
if $\alpha\in(1,2)$, we have
\begin{align*}
\|I_2\|_p&\leq4\|f\|_p\|a\|_\infty\int_{|y|>\eps}\nu(\dif y)+2\|\nabla f\|_p\|a\|_\infty\int_{|y|>\eps}|y|\nu(\dif y)\\
&\leq C_\eps\|f\|_p+C_\eps\|f\|^{\frac{1}{\alpha}}_{\alpha,p}\|f\|^{1-\frac{1}{\alpha}}_p
\leq \eps\|f\|_{\alpha,p}+C_\eps\|f\|_p.
\end{align*}
For $I_3(x)$, by Lemma \ref{Le1} we have
$$
\|I_3\|_p\leq C\gamma(\eps)\|(-\Delta)^{\alpha/2}f\|_p,
$$
where $\gamma(\eps)\to 0$ as $\eps\to 0$.

Now, combining the above calculations, we obtain the right hand side estimate in (\ref{EW1}). Moreover, we also have
$$
\|\cL^{a\nu}f\|_p\geq \|I_1\|_p-\|I_2\|_p-\|I_3\|_p\geq (a_0C_0-\eps-C\gamma(\eps))\|(-\Delta)^{\alpha/2}f\|_p-C_\eps\|f\|_p.
$$
Letting $\eps$ be small enough, we obtain the left hand side estimate in (\ref{EW1}).
\end{proof}

\section{Nonlocal linear parabolic equation}
In this section we fix a L\'evy measure $\nu$ satisfying {\bf (H$^{(\alpha)}_\nu$)}.
Let $\lambda: \mR_+\to\mR_+$ be a nonnegative and locally integrable function.
Let $N(\dif t,\dif x)$ be the Poisson random point measure with intensity measure
$\hat N(\dif t,\dif x):=\lambda(t)\dif t\nu(\dif x)$. Let
$\tilde N(\dif t,\dif x):=N(\dif t,\dif x)-\hat N(\dif t,\dif x)$ be the compensated random
martingale measure. Let $\vartheta:\mR_+\to\mR^d$ be a locally integrable function.
For $t\geq 0$, define
\begin{align}
X_t:=\int^t_0\vartheta(r)\dif r+\int^t_0\!\!\!\int_{B^{(\alpha)}}y\tilde N(\dif r,\dif y)
+\int^t_0\!\!\!\int_{\mR^d-B^{(\alpha)}}yN(\dif r,\dif y),\label{XX}
\end{align}
where $B^{(\alpha)}=\{x: |x|\leq 1\}$ if $\alpha=1$; $B^{(\alpha)}=\mR^d$ if $\alpha\in(1,2)$;
and $B^{(\alpha)}=\emptyset$ if $\alpha\in(0,1)$.

For $\varphi\in C^2_b(\mR^d)$, by It\^o's formula we have
\begin{align*}
&\mE\varphi(x+X_t-X_s)=\varphi(x)+\mE\int^t_s\vartheta(r)\cdot\nabla\varphi(x+X_r-X_s)\dif r\\
&\quad+\mE\int^t_s\!\!\!\int_{\mR^d}[\varphi(x+X_r-X_s+y)-\varphi(x+X_r-X_s)-y^{(\alpha)}\cdot\nabla\varphi(x+X_r-X_s)]\hat N(\dif r,\dif y).
\end{align*}
Thus, if we let
\begin{align}
\cT_{t,s}\varphi(x):=\cT^{\lambda\nu,\vartheta}_{t,s}\varphi(x):=\mE\varphi\left(x+X_t-X_s\right),\label{EY2}
\end{align}
then one sees that
$$
\p_t\cT_{t,s}\varphi=\cL^{\lambda(t)\nu}\cT_{t,s}\varphi+\vartheta(t)\cdot\nabla\cT_{t,s}\varphi.
$$

The following result is a slight extension of \cite[Theorem 4.2]{Zh1}.
\bt\label{Th1}
Assume {\bf (H$^{(\alpha)}_\nu$)} with $\alpha\in(0,2)$. Let $\vartheta:\mR_+\to\mR^d$ be a locally integrable
function and $\lambda:\mR_+\to[\lambda_0,\infty)$ be a measurable function, where $\lambda_0>0$. Let
$\cT^{\lambda\nu,\vartheta}_{t,s}$ be defined by (\ref{EY2}).
Then for any $p\in(1,\infty)$, there exists a constant $C=C(\lambda_0,\nu^{(\alpha)}_1,\nu^{(\alpha)}_2,\alpha,p,d)>0$
such that for any $T>0$ and $f\in L^p((0,T)\times\mR^d)$,
\begin{align}
\int^T_0\left\|\cL^{\nu}\int^t_0\cT^{\lambda\nu,\vartheta}_{t,s}f(s,\cdot)\dif s\right\|_p^p\dif t
\leq C\int^T_0\|f(t)\|^p_p\dif t.\label{Es44}
\end{align}
\et
\begin{proof}
Let $N^{(1)}(\dif t,\dif x)$ and $N^{(2)}(\dif t,\dif x)$ be two independent Poisson random point measures with
intensity measures $\hat N^{(1)}(\dif t,\dif x):=(\lambda(t)-\lambda_0)\dif t\nu(\dif x)$
and $\hat N^{(2)}(\dif t,\dif x):=\lambda_0\dif t\nu(\dif x)$ respectively.
Let $X^{(1)}_t$ be defined by (\ref{XX}) in terms of $N^{(1)}$, and $X^{(2)}_t$ be defined by
$$
X^{(2)}_t:=\int^t_0\!\!\!\int_{B^{(\alpha)}}y\tilde N^{(2)}(\dif r,\dif y)
+\int^t_0\!\!\!\int_{\mR^d-B^{(\alpha)}}yN^{(2)}(\dif r,\dif y).
$$
In fact, $X^{(2)}_t$ is the L\'evy process corresponding to the L\'evy measure $\lambda_0\nu(\dif x)$.
By It\^o's formula, we have
\begin{align}
\cT^{\lambda\nu,\vartheta}_{t,s}f(x)=\mE f(x+X^{(1)}_t-X^{(1)}_s+X^{(2)}_t-X^{(2)}_s)
=\mE\cT^{\lambda_0\nu,0}_{t,s}f(x+X^{(1)}_t-X^{(1)}_s).\label{EP1}
\end{align}
Thus, by \cite[Theorem 4.2]{Zh1}, we have
\begin{align*}
\int^T_0\left\|\cL^{\nu}\int^t_0\cT^{\lambda\nu,\vartheta}_{t,s}f(s,\cdot)\dif s\right\|_p^p\dif t
&\leq\mE\int^T_0\left\|\cL^{\nu}\int^t_0\cT^{\lambda_0\nu,0}_{t,s}f(s,\cdot+X^{(1)}_t-X^{(1)}_s)\dif s\right\|_p^p\dif t\\
&=\mE\int^T_0\left\|\cL^{\nu}\int^t_0\cT^{\lambda_0\nu,0}_{t,s}f(s,\cdot-X^{(1)}_s)\dif s\right\|_p^p\dif t\\
&\leq C\mE\int^T_0\left\|f(s,\cdot-X^{(1)}_s)\right\|_p^p\dif s
=C\int^T_0\|f(s)\|^p_p\dif s.
\end{align*}
The proof is complete.
\end{proof}

Consider the following time-dependent linear nonlocal parabolic system:
\begin{align}
\p_t u=\cL^{a(t)\nu}u+b^{(\alpha)}\cdot \nabla u+f,\ \ u(0)=\varphi,\label{EE5}
\end{align}
where $u,f:[0,1]\times\mR^d\to\mR^m$, $a:[0,1]\times\mR^d\times\mR^d\to\mR$ and
$b:[0,1]\times\mR^d\to\mR^d$ are Borel measurable functions, and
\begin{align}
b^{(\alpha)}(t,x)=1_{\alpha\in[1,2)}b(t,x).\label{BB}
\end{align}
We make the following assumptions on $a$ and $b$:
\begin{enumerate}[{\bf (H$^a_\nu$)}]
\item For each $t\geq 0$, $a(t)\in\sA_1$ satisfies
$$
\sup_{t\in[0,1]}\|a(t)\|_{\sA_1}<+\infty,\ \ a_0\leq a(t,x,0)\leq a_1,
$$
where $a_0,a_1>0$, and for all $0<r<R<+\infty$,
\begin{align}
1_{\alpha=1}\int_{r\leq |y|\leq R}ya(t,x,y)\nu(\dif y)=0.\label{EL222}
\end{align}
\end{enumerate}
\begin{enumerate}[{\bf (H$^b$)}]
\item For all $t\geq 0$ and $x,y\in\mR^d$,
$$
|b^{(\alpha)}(t,x)-b^{(\alpha)}(t,y)|\leq 1_{\alpha=1}\omega_b(|x-y|)+1_{\alpha\in(1,2)}C_b,
$$
where $\omega_b:\mR^+\to\mR^+$ is an increasing function with $\lim_{s\downarrow 0}\omega_b(s)=0$.
\end{enumerate}

Let us first prove the following apriori estimate by the method of freezing the coefficients (cf. \cite[Lemma 5.1]{Zh1}).
\bl\label{Le7}
Suppose that $a(t,x,y)=a(t,x)$ is independent of $y$  and satisfies {\bf (H$^a_\nu$)},
and $b$ satisfies {\bf (H$^b$)}.
Let $p>1$ and not equal to $\frac{\alpha}{\alpha-1}$ if $\alpha\in(1,2)$, and let $f\in \mY^{0,p}$ and
$u\in\mX^{\alpha,p}$ satisfy (\ref{EE5}). Then for all $t\in[0,1]$,
\begin{align}
\|u\|_{\mX^{\alpha,p}_t}\leq C\left(\|u(0)\|_{\mW^{\alpha-\frac{\alpha}{p},p}}
+\|f\|_{\mY^{0,p}_t}\right),\label{Es888}
\end{align}
where $C$ depends only on $a_0,a_1,\|a\|_{\sA_1}, \|b\|_\infty,d,p,\alpha$ and $\omega_b$.
\el
\begin{proof}
Let $(\rho_\eps)_{\eps\in(0,1)}$ be a family of mollifiers in $\mR^d$, i.e., $\rho_\eps(x)=\eps^{-d}\rho(\eps^{-1}x)$,
where $\rho\in C^\infty_0(\mR^d)$ with $\int\rho=1$ is nonnegative. Define
$$
u_\eps(t):=u(t)*\rho_\eps,\ \ a_\eps(t):=a(t)*\rho_\eps,\ \ b_\eps(t):=b(t)*\rho_\eps,\ \ f_\eps(t):=f(t)*\rho_\eps,
$$
where $*$ stands for the convolution.
Taking convolutions for both sides of (\ref{EE5}), we have
\begin{align}
\p_tu_\eps=\cL^{a_\eps(t)\nu}u_\eps+b^{(\alpha)}_\eps\cdot\nabla u_\eps+h_\eps,\label{Eq333}
\end{align}
where
$$
h_\eps:=f_\eps+(\cL^{a(t)\nu} u)*\rho_\eps-\cL^{a_\eps(t)\nu} u_\eps
+(b^{(\alpha)}\cdot\nabla u)*\rho_\eps-b^{(\alpha)}_\eps\cdot\nabla u_\eps.
$$
By the assumption, it is easy to see that for all $\eps\in(0,1)$ and $t\in[0,1]$ and $x,y\in\mR^d$,
\begin{align}
|a_\eps(t,x)-a_\eps(t,y)|\leq\omega^{(1)}_a(|x-y|),\ \
|b_\eps(t,x)-b_\eps(t,y)|\leq 1_{\alpha=1}\omega_b(|x-y|)+1_{\alpha\in(1,2)}C_b,\label{Eq5}
\end{align}
and
$$
|a_\eps(t,x)-a(t,x)|\leq\omega^{(1)}_a(\eps),\ \ |b_\eps(t,x)-b(t,x)|\leq
1_{\alpha=1}\omega_b(\eps)+1_{\alpha\in(1,2)}C_b.
$$
Moreover, by the property of convolutions, we also have
$$
\lim_{\eps\to0}\int^1_0\|h_\eps(t)-f(t)\|^p_p\dif t=0.
$$
Below, we use the method of freezing the coefficients to prove that for all $t\in[0,1]$,
\begin{align}
\|u_\eps\|_{\mX^{\alpha,p}_t}\leq C_{t,p}\left(\|u_\eps(0)\|_{\mW^{\alpha-\frac{\alpha}{p}},p}
+\|h_\eps\|_{\mY^{0,p}_t}\right),\label{Ep7}
\end{align}
where the constant $C$ is independent of $\eps$. After proving this estimate, (\ref{Es888})
immediately follows by taking limits for (\ref{Ep7}).

For simplicity of notation, we drop the subscript $\eps$ below.
Fix $\delta>0$ being small enough, whose value will be determined below.
Let $\zeta$ be a smooth function with support in $B_\delta$ and $\|\zeta\|_p=1$. For $z\in\mR^d$, set
$$
\zeta_z(x):=\zeta(x-z),\ \ \lambda^a_z(t):=a(t,z),\ \ \vartheta^b_z(t):=1_{\alpha=1}b(t,z).
$$
Multiplying both sides of (\ref{Eq333}) by $\zeta_z$, we have
$$
\p_t(u\zeta_z)=\lambda^a_z\cL^\nu (u\zeta_z)+\vartheta^b_z\cdot\nabla(u\zeta_z)+g^\zeta_z,
$$
where
$$
g^\zeta_z:=(a-\lambda^a_z)\cL^\nu u\zeta_z
+\lambda^a_z(\cL^\nu u\zeta_z-\cL^\nu(u\zeta_z))
+(b^{(\alpha)}-\vartheta^b_z)\cdot\nabla (u\zeta_z)-ub^{(\alpha)}\cdot\nabla\zeta_z+h\zeta_z.
$$
Let $\cT^{\lambda^a_z\nu,\vartheta^b_z}_{t,s}$ be defined by (\ref{EY2}) in terms of $\lambda^a_z\nu$ and $\vartheta^b_z$.
By Duhamel's formula, $u\zeta_z$ can be written as
$$
u\zeta_z(t,x)=\cT^{\lambda^a_z\nu,\vartheta^b_z}_{t,0}(u(0)\zeta_z)(x)+\int^t_0\cT^{\lambda^a_z\nu,\vartheta^b_z}_{t,s} g^\zeta_z(s,x)\dif s,
$$
and so that for any $T\in[0,1]$,
\begin{align*}
\int^T_0\|\cL^\nu(u\zeta_z)(t)\|^p_p\dif t&\leq
2^{p-1}\left(\int^T_0\|\cL^\nu\cT^{\lambda^a_z\nu,\vartheta^b_z}_{t,0}(u(0)\zeta_z)\|_p^p\dif t
+\int^T_0\left\|\cL^\nu\int^t_0\cT^{\lambda^a_z\nu,\vartheta^b_z}_{t,s}
g^\zeta_z(s)\dif s\right\|^p_p\dif t\right)\\
&=:2^{p-1}(I_1(T,z)+I_2(T,z)).
\end{align*}
For $I_1(T,z)$, by (\ref{EP1}) and $\|\cL^\nu f(\cdot+z)\|_p=\|\cL^\nu f\|_p$, we have
\begin{align}
\int^T_0\|\cL^\nu\cT^{\lambda^a_z\nu,\vartheta^b_z}_{t,0}(u(0)\zeta_z)\|^p_p\dif t
=\int^T_0\left\|\cL^\nu\cT^{a_0\nu,0}_{t,0}(u(0)\zeta_z)\right\|^p_p\dif t
\leq C\|u(0)\zeta_z\|^p_{\mW^{\alpha-\frac{\alpha}{p},p}},\label{ET2}
\end{align}
where the last step is due to \cite[p.96 Theorem 1.14.5]{Tr} and \cite[Corollary 4.5]{Zh1}.
Thus, by definition (\ref{Eq404}), it is easy to see that
$$
\int_{\mR^d}I_1(T,z)\dif z\leq C\int_{\mR^d}\|u(0)\zeta_z\|^p_{\mW^{\alpha-\frac{\alpha}{p},p}}\dif z
\leq C\Big(\|u(0)\|^p_{\mW^{\alpha-\frac{\alpha}{p},p}}\|\zeta\|_p^p
+\|u(0)\|_p^p\|\zeta\|^p_{\mW^{\alpha-\frac{\alpha}{p},p}}\Big).
$$
For $I_2(T,z)$, by Theorem \ref{Th1}, we have
\begin{align*}
I_2(T,z)&\leq C\int^T_0\|g^\zeta_z(s)\|^p_p\dif s\leq
C\int^T_0\|((a-\lambda^a_z)(\cL^\nu u\zeta_z))(s)\|^p_p\dif s\\
&\qquad+C\int^T_0\|\lambda^a_z(\cL^\nu(u\zeta_z)-\cL^\nu u\zeta_z)(s)\|^p_p\dif s\\
&\qquad+C\int^T_0\|((b^{(\alpha)}-\vartheta^b_z)\cdot\nabla(u\zeta_z))(s)\|^p_p\dif s\\
&\qquad+C\int^T_0\|(ub^{(\alpha)}\cdot\nabla\zeta_z)(s)\|^p_p\dif s+C\int^T_0\|h\zeta_z(s)\|^p_p\dif s\\
&=:I_{21}(T,z)+I_{22}(T,z)+I_{23}(T,z)+I_{24}(T,z)+I_{25}(T,z).
\end{align*}
For $I_{21}(T,z)$, by (\ref{Eq5}) and $\|\zeta\|_p=1$, we have
$$
\int_{\mR^d}I_{21}(T,z)\dif z\leq C\omega^{(1)}_a(\delta)^p
\int_{\mR^d}\!\!\int^T_0\|(\cL^\nu u\zeta_z)(s)\|^p_p\dif s\dif z
= C\omega^{(1)}_a(\delta)^p\int^T_0\|\cL^\nu u(s)\|^p_p\dif s.
$$
For $I_{22}(T,z)$, using (\ref{Le2}) and as in the proof of \cite[Lemma 2.5]{Zh1},
for any $\beta\in(0\vee(\alpha-1),\alpha)$, we have
\begin{align*}
\int_{\mR^d}I_{22}(T,z)\dif z&\leq
Ca_1\int^T_0\!\!\!\int_{\mR^d}\|(\cL^\nu(u\zeta_z)-\cL^\nu u\zeta_z)(s)\|^p_p\dif z\dif s\\
&\leq C\int^T_0\|u(s)\|^p_p\dif s+C\int^T_0\|(-\Delta)^{\beta/2}u(s)\|^{p}_p\dif s\\
&\leq C\int^T_0\|u(s)\|^p_p\dif s+\frac{1}{4^p}\int^T_0\|\cL^\nu u(s)\|^p_p\dif s,
\end{align*}
where the last step is due to the interpolation inequality, Young's inequalities and Theorem \ref{Th3}.

For $I_{23}(T,z)$, as above we have
$$
\int_{\mR^d}I_{23}(T,z)\dif z\leq C1_{\alpha=1}\omega_b(\delta)^p\left(\int^T_0\|\nabla u(s)\|^p_p\dif s+\|\nabla\zeta\|^p_p\int^T_0\|u(s)\|^p_p\dif s\right).
$$
Moreover, it is easy to see that
\begin{align*}
\int_{\mR^d}I_{24}(T,z)\dif z&\leq C\|b\|^p_\infty\|\nabla\zeta\|_p^p\int^T_0\|u(s)\|^p_p\dif s,\\
\int_{\mR^d}I_{25}(T,z)\dif z&\leq C\int^T_0\|h(s)\|^p_p\dif s.
\end{align*}
Combining the above calculations, we get
\begin{align*}
\int^T_0\|\cL^\nu u(s)\|^p_p\dif s&=\int^T_0\!\!\!\int_{\mR^d}\|\cL^\nu u(s)\cdot\zeta_z\|^p_p\dif z\dif s
\leq 2^{p-1}\int^T_0\!\!\!\int_{\mR^d}\|\cL^\nu(u\zeta_z)(s)\|^p_p\dif z\dif s\\
&\quad+2^{p-1}\int^T_0\!\!\!\int_{\mR^d}\|(\cL^\nu u\zeta_z-\cL^\nu(u\zeta_z))(s)\|^p_p\dif z\dif s\\
&\leq C\|u(0)\|^p_{\mW^{\alpha-\frac{\alpha}{p},p}}+\Big(\frac{1}{4}
+C(\omega^{(1)}_a(\delta)^p+\omega_b(\delta)^p)\Big)\int^T_0\|\cL^\nu u(s)\|^p_p\dif s\\
&\quad+C\int^T_0\|u(s)\|^p_p\dif s+C\int^T_0\|h(s)\|^p_p\dif s.
\end{align*}
Choosing $\delta_0>0$ being small enough so that
$$
C(\omega^{(1)}_a(\delta_0)^p+\omega_b(\delta_0)^p)\leq \frac{1}{4},
$$
we obtain that for all $T\in[0,1]$,
\begin{align}
\int^T_0\|\cL^\nu u(s)\|^p_p\dif s\leq C\|u(0)\|^p_{\mW^{\alpha-\frac{\alpha}{p},p}}+
C\int^T_0\|u(s)\|^p_p\dif s+C\int^T_0\|h(s)\|^p_p\dif s.\label{EL1}
\end{align}
On the other hand, by (\ref{Eq333}), it is easy to see that for all $t\in[0,1]$,
$$
\|u(t)\|^p_p\leq C\|u(0)\|_p^p+C1_{\alpha\in[1,2)}\int^t_0 \|\nabla u(s)\|^p_p\dif s+C\int^t_0 \|h(s)\|^p_p\dif s,
$$
which together with (\ref{EL1}) and Gronwall's inequality yields that for all $t\in[0,1]$,
\begin{align}
\sup_{s\in[0,t]}\|u(s)\|^p_p+\int^t_0\|\cL^\nu u(s)\|^p_p\dif s\leq C
\left(\|u(0)\|^p_{\mW^{\alpha-\frac{\alpha}{p},p}}+\int^t_0 \|h(s)\|^p_p\dif s\right).\label{EI1}
\end{align}
From equation (\ref{Eq333}), we also have
\begin{align*}
\int^t_0\|\p_su(s)\|^p_p\dif s&\leq C\left(\|a\|_\infty^p\int^t_0\|\cL^\nu u(s)\|^p_p\dif s
+\|b^{(\alpha)}\|_\infty^p\int^t_0\|\nabla u(s)\|^p_p\dif s
+\int^t_0\|h(s)\|^p_p\dif s\right),
\end{align*}
which together with (\ref{EI1}) and (\ref{EW88}) gives (\ref{Ep7}), and therefore (\ref{Es888}).
\end{proof}

We now prove the following main result of this paper.
\bt
Suppose {\bf (H$^{(\alpha)}_\nu$)}, {\bf (H$^a_\nu$)} and {\bf (H$^b$)} and for some $k\in\mN\cup\{0\}$,
$$
|\nabla^j_xa(t,x,y)|+|\nabla^j_xb(t,x)|\leq C_j, \  \ j=0,\cdots,k.
$$
For given $p\in(1,\infty)$ not equal to $\frac{\alpha}{\alpha-1}$ when $\alpha\in(1,2)$,
and for $\varphi\in\mW^{k+\alpha-\frac{\alpha}{p},p}$, there exists
a unique $u\in\mX^{k+\alpha,p}$ satisfying equation (\ref{EE5}). Moreover, for all $t\in[0,1]$,
\begin{align}
\|u\|_{\mX^{k+\alpha,p}_t}\leq C_{k,p}\left(\|\varphi\|_{\mW^{k+\alpha-\frac{\alpha}{p},p}}
+\|f\|_{\mY^{k,p}_t}\right),\label{Es8}
\end{align}
where $C_{0,p}$ depends only on $a_0,a_1,\|a\|_{\sA_1},\|b\|_\infty,d,p,\alpha$ and $\omega_b$.
\et
\begin{proof}
The strategy is to prove the apriori estimate (\ref{Es8}) and then use the continuity method (cf. \cite{Kr3, Zh3}).

(Step 1) Let us first rewrite equation (\ref{EE5}) as
$$
\p_tu(t,x)=a(t,x,0)\cL^{\nu} u(t,x)+b^{(\alpha)}(t,x)\cdot\nabla u(t,x)+\tilde f(t,x),\label{EE55}
$$
where
$$
\tilde f(t,x):=f(t,x)+\int_{\mR^d}\cJ^{(\alpha)}_{u(t,\cdot)}(x,y)(a(t,x,y)-a(t,x,0))\nu(\dif y)
$$
and
$$
\cJ^{(\alpha)}_{u(t,\cdot)}(x,y):=u(t,x+y)-u(t,x)-y^{(\alpha)}\cdot\nabla u(t,x).
$$
Notice that by Lemma \ref{Le1},
\begin{align*}
\|\tilde f(t)\|_p&\leq\|f(t)\|_p+
\left\|\int_{|y|>\eps}\cJ^{(\alpha)}_{u(t,\cdot)}(\cdot,y)(a(t,\cdot,y)-a(t,\cdot,0))\nu(\dif y)\right\|_p\\
&\quad+\left\|\int_{|y|\leq \eps}\cJ^{(\alpha)}_{u(t,\cdot)}(\cdot,y)(a(t,\cdot,y)-a(t,\cdot,0))\nu(\dif y)\right\|_p\\
&\leq\|f(t)\|_p+2a_1\left(\|u(t)\|_p\nu(B^c_\eps)+1_{\alpha\in(1,2)}\|\nabla u(t)\|_p\right)\\
&\quad+C_\eps\|u(t)\|_p+C\gamma_0(\eps)\|(-\Delta)^{\frac{\alpha}{2}}u(t)\|_p\\
&\leq\|f(t)\|_p+C_\eps\|u(t)\|_p+\gamma_1(\eps)\|(-\Delta)^{\frac{\alpha}{2}}u(t)\|_p,
\end{align*}
where the last step is due to the interpolation inequality and Young's inequalities, and
$$
\gamma_0(\eps):=\int^\eps_0\frac{\omega^{(0)}_a(r)}{r}\dif r,\ \ \gamma_1(\eps):=\eps+C\gamma_0(\eps).
$$
By Lemma \ref{Le7}, we have
$$
\|u\|_{\mX^{\alpha,p}_t}\leq C_1\|\varphi\|_{\mW^{\alpha-\frac{\alpha}{p},p}}+C_2\|\tilde f\|_{\mY^{0,p}_t}.
$$
In particular, for all $t\in[0,1]$,
\begin{align*}
\sup_{s\in[0,t]}\|u(s)\|_p+\int^t_0\|(-\Delta)^{\frac{\alpha}{2}}u(s)\|_p^p\dif s
&\leq C_1\|\varphi\|_{\mW^{\alpha-\frac{\alpha}{p},p}}
+\gamma_1   (\eps)\int^t_0\|(-\Delta)^{\frac{\alpha}{2}}u(s)\|_p^p\dif s\\
&\quad+C_2\int^t_0\|u(s)\|^p_p\dif s+C_2\int^t_0\|f(s)\|^p_p\dif s.
\end{align*}
Letting $\eps$ be small enough and using Gronwall's inequality, we obtain (\ref{Es8}) with $k=0$.

(Step 2) We now estimate the higher order derivatives. Write
$$
w^{(n)}(t,x):=\nabla^n u(t,x).
$$
By the chain rule, one can see that
$$
\p_tw^{(n)}=\cL^{a\nu}w^{(n)}+b^{(\alpha)}\cdot\nabla w^{(n)}+g^{(n)},
$$
where
$$
g^{(n)}:=\nabla^n f+\sum_{j=1}^n\frac{n!}{(n-j)! j!}\left(\cL^{(\nabla^j_x a)\nu}(\nabla^{n-j}u)
+\nabla^j b^{(\alpha)}\cdot\nabla^{n-j+1}u\right)
$$
and
$$
\cL^{(\nabla^j_x a)\nu}(\nabla^{n-j}u)(t,x)=\int_{\mR^d}\cJ^{(\alpha)}_{\nabla^{n-j}u(t,\cdot)}(x,y)
\nabla^j_xa(t,x,y)\nu(\dif y).
$$
By Step 1, we know that
\begin{align}
\|w^{(n)}\|_{\mX^{\alpha,p}_t}\leq C\left(\|w^{(n)}(0)\|_{\mW^{\alpha-\frac{\alpha}{p},p}}
+\|g^{(n)}\|_{\mY^{0,p}_t}\right).\label{EU1}
\end{align}
By Minkowskii's inequality, we have
\begin{align*}
\|\cL^{(\nabla^j_x a)\nu}(\nabla^{n-j}u)(t)\|_p
&\leq C_j\int_{\mR^d}\|\nabla^{n-j}u(t,\cdot+y)-\nabla^{n-j}u(t,\cdot)
-y^{(\alpha)}\cdot\nabla\nabla^{n-j}u(t,\cdot)\|_p\nu(\dif y)\\
&\leq C_j\int_{|y|\geq1}\Big(2\|\nabla^{n-j}u(t)\|_p+1_{\alpha\in(1,2)}|y|\cdot\|\nabla^{n-j+1}u(t)\|_p\Big)\nu(\dif y)\\
&\quad+C_j1_{\alpha\in(0,1)}\int_{|y|\leq 1}\|\nabla^{n-j}u(t,\cdot+y)-\nabla^{n-j}u(t,\cdot)\|_p\nu(\dif y)\\
&\quad+C_j1_{\alpha\in[1,2)}\int_{|y|\leq 1}\|\nabla^{n-j}u(t,\cdot+y)-\nabla^{n-j}u(t,\cdot)
-y\cdot\nabla\nabla^{n-j}u(t,\cdot)\|_p\nu(\dif y)\\
&\leq2 C_j\nu(B^c_1)\|\nabla^{n-j}u(t)\|_p+C_j\|\nabla^{n-j+1}u(t)\|_p\\
&\qquad\qquad\qquad\times\left(\int_{B^c_1}|y|1_{\alpha\in(1,2)}\nu(\dif y)
+\int_{B_1}|y|1_{\alpha\in(0,1)}\nu(\dif y)\right)\\
&\quad+C_j1_{\alpha\in[1,2)}\int_{|y|\leq 1}|y|\left(\int^1_0\|\nabla^{n-j+1}u(t,\cdot+sy)
-\nabla^{n-j+1}u(t,\cdot)\|_p\dif s\right)\nu(\dif y)\\
&\leq C\|\nabla^{n-j}u(t)\|_p+C\|\nabla^{n-j+1}u(t)\|_p
+C\|\nabla^{n-j+1}u(t)\|_{\beta,p} \int_{B_1}|y|^{1+\beta}1_{\alpha\in[1,2)}\nu(\dif y),
\end{align*}
where $\beta\in((\alpha-1)\vee 0,1)$ and the last step is due to (\ref{Le2}).

Hence, by the assumptions, we obtain
\begin{align*}
\|g^{(n)}\|^p_{\mY^{0,p}_t}\leq\|f\|^p_{\mY^{n,p}_t}+C\|u\|^p_{\mY^{n,p}_t}
+C_t\|u\|^p_{\mY^{n+\beta,p}_t}1_{\alpha\in[1,2)}.
\end{align*}
Summing over $n$ from $0$ to $k$ for (\ref{EU1}) yields
\begin{align*}
\|u(t)\|^p_{\mW^{k,p}}+\int^t_0\|u(s)\|^p_{\mW^{k+\alpha,p}}\dif s&\leq
C\|\varphi\|_{\mW^{k+\alpha-\frac{\alpha}{p},p}}+C1_{\alpha\in[1,2)}\int^t_0\|u(s)\|^p_{\mW^{k+\beta,p}}\dif s\\
&\quad+C\int^t_0\|f(s)\|^p_{\mW^{k,p}}\dif s+C\int^t_0\|u(s)\|^p_{\mW^{k,p}}\dif s\\
&\leq C\|\varphi\|_{\mW^{k+\alpha-\frac{\alpha}{p},p}}+C1_{\alpha\in[1,2)}\int^t_0
\|u(s)\|^{p\beta /\alpha}_{\mW^{k+\alpha,p}}\|u(s)\|^{p(1-\beta/\alpha)}_{\mW^{k,p}}\dif s\\
&\quad+C\int^t_0\|f(s)\|^p_{\mW^{k,p}}\dif s+C\int^t_0\|u(s)\|^p_{\mW^{k,p}}\dif s\\
&\leq C\|\varphi\|_{\mW^{k+\alpha-\frac{\alpha}{p},p}}+\frac{1}{2}1_{\alpha\in[1,2)}\int^t_0
\|u(s)\|^p_{\mW^{k+\alpha,p}}\dif s\\
&\quad+C\int^t_0\|f(s)\|^p_{\mW^{k,p}}\dif s+C\int^t_0\|u(s)\|^p_{\mW^{k,p}}\dif s,
\end{align*}
which then gives (\ref{Es8}) by Gronwall's inequality.

(Step 3) For $\lambda\in[0,1]$, define an operator
$$
U_\lambda:=\p_t-\lambda \cL^{a\nu}-\lambda b^{(\alpha)}\cdot\nabla-(1-\lambda)\cL^\nu.
$$
By (\ref{EW1}), it is easy to see that
\begin{align}
U_\lambda: \mX^{k+\alpha,p}\to \mY^{k,p}.\label{Eq6}
\end{align}
For given $\varphi\in\mW^{k+\alpha-\frac{\alpha}{p},p}$, let $\mX^{k+\alpha,p}_\varphi$ be the space
of all functions $u\in\mX^{k+\alpha,p}$ with $u(0)=\varphi$. It is clear that $\mX^{k+\alpha,p}_\varphi$
is a complete metric space with respect to the metric $\|\cdot\|_{\mX^{k+\alpha,p}}$.
For $\lambda=0$ and $f\in \mY^{k,p}$, it is well-known that
there is a unique $u\in\mX^{k+\alpha,p}_\varphi$ such that
$$
U_0 u=\p_t u-\cL^\nu u=f.
$$
In fact, by Duhamel's formula, the unique solution can be represented by
$$
u(t,x)=\cT^{\nu,0}_{t,0}\varphi(x)+\int^t_0\cT^{\nu,0}_{t,s}f(s,x)\dif s,
$$
where $\cT^{\nu,0}_{t,s}$ is defined by (\ref{EY2}).
Suppose now that for some $\lambda_0\in[0,1)$, and for any $f\in \mY^{k,p}$, the equation
$$
U_{\lambda_0}u=f
$$
admits a unique solution $u\in\mX^{k+\alpha,p}_\varphi$.
Thus, for fixed $f\in\mY^{k,p}$ and $\lambda\in[\lambda_0,1]$, and for any $u\in\mX^{k+\alpha,p}_\varphi$,
by (\ref{Eq6}), the equation
\begin{align}
U_{\lambda_0}w=f+(U_{\lambda_0}-U_\lambda)u\label{Eq0}
\end{align}
admits a unique solution $w\in\mX^{k+\alpha,p}_\varphi$.
Introduce an operator
$$
w=Q_\lambda u.
$$
We now use the apriori estimate (\ref{Es8}) to show that there exists an $\eps>0$
independent of $\lambda_0$ such that for all $\lambda\in[\lambda_0,\lambda_0+\eps]$,
$$
Q_\lambda:\mX^{k+\alpha,p}_\varphi\to\mX^{k+\alpha,p}_\varphi
$$
is a contraction operator.

Let $u_1,u_2\in\mX^{k+\alpha,p}_\varphi$ and $w_i=Q_{\lambda}u_i,i=1,2$. By equation (\ref{Eq0}), we have
\begin{align*}
U_{\lambda_0}(w_1-w_2)=(U_{\lambda_0}-U_\lambda)(u_1-u_2)=(\lambda_0-\lambda)(\cL^{(a-1)\nu}+b^{(\alpha)}\cdot\nabla)(u_1-u_2).
\end{align*}
By (\ref{Es8}) and (\ref{EW1}), it is not hard to see that
\begin{align*}
\|Q_\lambda u_1-Q_\lambda u_2\|_{\mX^{k+\alpha,p}}
&\leq C_{k,p}|\lambda_0-\lambda|\cdot\|(\cL^{(a-1)\nu}+b^{(\alpha)}\cdot\nabla)(u_1-u_2)\|_{\mY^{k,p}}\\
&\leq C_0|\lambda_0-\lambda|\cdot\|u_1-u_2\|_{\mX^{k+\alpha,p}},
\end{align*}
where $C_0$ is independent of $\lambda,\lambda_0$ and $u_1,u_2$. Taking $\eps=1/(2C_0)$, one sees that
$$
Q_\lambda: \mX^{k+\alpha,p}_\varphi\to \mX^{k+\alpha,p}_\varphi
$$
is a $1/2$-contraction operator. By the fixed point theorem, for each $\lambda\in[\lambda_0,\lambda_0+\eps]$,
there exists a unique $u\in\mX^{k+\alpha,p}_\varphi$ such that
$$
Q_\lambda u=u,
$$
which means that
$$
U_\lambda u=f.
$$
Now starting from $\lambda=0$, after repeating the above construction
$[\frac{1}{\eps}]+1$-steps, one obtains that for any $f\in\mY^{k,p}$,
$$
U_1u=f
$$
admits a unique solution $u\in\mX^{k+\alpha,p}_\varphi$.
\end{proof}

\end{document}